# An intrusion detection system in internet of things using grasshopper optimization algorithm and machine learning algorithms

Shiva Sattarpour[1], Ali Barati[1*], Hamid Barati[1]

1 Department of Computer Engineering, Dezful Branch, Islamic Azad University, Dezful, Iran

**Abstract**

The Internet of Things (IoT) has emerged as a foundational paradigm supporting a range of applications, including healthcare, education, agriculture, smart homes, and, more recently, enterprise systems. However, significant advancements in IoT networks have been impeded by security vulnerabilities and threats that, if left unaddressed, could hinder the deployment and operation of IoT-based systems. Detecting unwanted activities within the IoT is crucial, as it directly impacts confidentiality, integrity, and availability. Consequently, intrusion detection has become a fundamental research area and the focus of numerous studies. An intrusion detection system (IDS) is essential to the IoT's alarm mechanisms, enabling effective security management. This paper examines IoT security and introduces an intelligent two-layer intrusion detection system for IoT. Machine learning techniques power the system's intelligence, with a two-layer structure enhancing intrusion detection. By selecting essential features, the system maintains detection accuracy while minimizing processing overhead. The proposed method for intrusion detection in IoT is implemented in two phases. In the first phase, the Grasshopper Optimization Algorithm (GOA) is applied for feature selection. In the second phase, the Support Vector Machine (SVM) algorithm is used to detect intrusions. The method was implemented in MATLAB, and the NSL-KDD dataset was used for evaluation. Simulation results show that the proposed method improves accuracy compared to other approaches.

## 1. Introduction

The rapid evolution of various technologies, such as sensors, automatic detection and tracking, embedded computing, wireless communication, broadband Internet access, and distributed services, has significantly enhanced the potential to integrate smart objects into our daily lives through the Internet [1]. This convergence of the Internet and intelligent, interconnected objects defines what we know as the Internet of Things (IoT). IoT holds promise across numerous application areas, including logistics, industrial processes, public safety, home automation, environmental monitoring, and healthcare [2]. However, the integration of real-world objects with the Internet introduces significant cybersecurity challenges that can impact various aspects of daily life [3].

To address these challenges, ongoing research and projects are focused on enhancing IoT security through methods that ensure data confidentiality and authentication, manage access control within IoT networks, and promote privacy and trust among users and devices [4]. These efforts include enforcing security and privacy policies to mitigate potential risks [5]. Nevertheless, IoT networks remain vulnerable to a range of sophisticated attacks intended to disrupt or compromise these networks [6]. Therefore, beyond preventive measures, it is essential to implement defensive strategies that detect potential attackers. Intrusion Detection Systems (IDS) are fundamental to this approach, aiming to identify and counter unauthorized attempts to access or disrupt information systems [7].

An IDS serves as a crucial tool in protecting both conventional networks and IoT-based systems [8]. It continuously monitors the activities within a host or network, promptly notifying the management system in the event of any security breaches [9]. Effective intrusion detection not only guards against unauthorized access but also strengthens overall network resilience by quickly identifying suspicious activities [10].

In this paper, we introduce a novel approach for detecting attackers in IoT environments by employing dynamic, behavior-based techniques integrated with data mining methods. The proposed framework operates in two primary phases: initially, the Grasshopper Optimization Algorithm (GOA) is used for feature reduction, enabling the system to identify and retain the most relevant features for intrusion detection. In the subsequent phase, a Support Vector Machine (SVM) classifier is applied to categorize the data samples, facilitating precise detection of potential intrusions. This dual-phase approach harnesses the strengths of both feature optimization and machine learning, significantly enhancing the accuracy and dependability of intrusion detection in IoT networks. By reducing irrelevant features and focusing on critical ones, this method not only improves detection accuracy but also optimizes processing efficiency, making it particularly suitable for the complex and data-intensive nature of IoT environments.

Here are the main innovations and contributions of the proposed method in this article:

- Hybrid Detection Approach: The method combines the Grasshopper Optimization Algorithm (GOA) for feature selection with a Support Vector Machine (SVM) classifier, enhancing accuracy by selecting only the most relevant features, which optimizes resource usage in IoT environments.

- Efficient Feature Selection: By implementing GOA for feature reduction, the model efficiently identifies significant attributes from the NSL-KDD dataset, minimizing computational requirements and maintaining high detection precision.

- Enhanced Intrusion Detection Performance: The dual-phase approach with GOA and SVM significantly improves intrusion detection rates, achieving high True Positive Rates (TPR) while maintaining low False Positive Rates (FPR).

- Adaptability for IoT Constraints: The proposed model addresses the resource limitations typical of IoT networks, ensuring low processing overhead and accurate threat detection suited to the heterogeneous and dynamic nature of IoT devices and environments.
- Comprehensive Evaluation: The method demonstrates superior performance compared to existing approaches, achieving higher accuracy, TPR, and lower FPR, as validated through simulations on the NSL-KDD dataset.

The remainder of this paper is organized as follows: In Sect. 2, the previous works are reviewed, and the advantages and disadvantages of each of them are also stated. The problem statement is described in Sect. 3. The proposed method is described in Sect. 4. The proposed method is compared with other schemes, and the results are displayed in the form of tables and charts in Sect. 5. Finally, conclusions are presented.

## 2. Related Works

In recent research on intrusion detection in Internet of Things (IoT) networks, various methods have been proposed to identify abnormal behaviours and enhance security across connected devices [11]. Given the unique architecture and extensive data flows within IoT environments, intrusion detection systems (IDS) for IoT must address specific challenges, such as resource constraints, heterogeneous devices, and a dynamic network topology [12]. These methods leverage a combination of behaviour-based techniques, machine learning algorithms, and data mining processes to detect suspicious activities and prevent unauthorized access [13]. By focusing on feature optimization and employing intelligent classification algorithms, these systems can enhance the accuracy and efficiency of intrusion detection while minimizing resource consumption, making them suitable for the diverse and large-scale nature of IoT networks [14].

Rajeshkumar et al. [15], developed a trust-based, energy-aware intrusion detection system that utilizes a Kalman filter to predict the trustworthiness of nodes and employs a double-ACK algorithm to reduce control costs and increase path reliability. This method successfully improved the packet delivery ratio and reduced network delay compared to previous approaches. However, it incurred higher control costs.

Qu et al. [16], introduced an intrusion detection system that combines a support vector machine (SVM) with a mean shift clustering algorithm for detecting intrusions and abnormal behaviors. This method achieved a high detection rate and low error rate, along with the capability to update features. However, it imposed additional processing overhead on the network.

Liu et al. [17], proposed a system aimed at intrusion detection in complex, large-scale networks, leveraging an improved social spider optimization algorithm and distributed fuzzy clustering. While this method demonstrated high efficiency, its limitations include the lack of deployment in real-world networks and insufficient focus on the detection accuracy for smaller sample sizes.

Pardhan et al. [18] introduced LB-IDS, a trust-based intrusion detection system enhancing WSN security by reducing message complexity, memory, and energy usage. LB-IDS detects transmembrane, reverse manipulation, hole-hole, and parasitic intrusions across multiple layers. Despite improved accuracy, LB-IDS lacks real-world testing, especially with outdoor transceiver modules.

Cauteruccio et al. [19] developed an anomaly detection method for WSNs using machine learning and multi-parameter editing distance, relying on temperature, humidity, and light sensors. The system combines short- and long-term strategies, achieving accurate and timely detection while adapting to environmental changes. Future improvements include enhancing automation and optimizing thresholds.

Amouri et al. [20] introduced AMoF, a machine learning-based intrusion detection system for mobile WSNs with a cross-layer architecture. AMoF operates in two stages: DS nodes filter data to create CCI reports, and SN nodes analyze this data via linear regression. Using a sliding window algorithm, AMoF improves attack detection accuracy.

Ferrag et al. [21] introduced RDTIDS, an intrusion detection system combining decision tree and rule-based methods (REP Tree, JRip, Forest PA) for WSNs. RDTIDS uses a hierarchical model for binary and multilevel classification. Through sequential training of three classifiers, it improves recognition accuracy and reduces data imbalance and computational errors.

Safaldin et al. [22] introduced an improved binary gray wolf optimization (IBGWO) algorithm with a support vector machine (SVM) for intrusion detection in WSNs. This approach enhances detection accuracy by optimizing the algorithm's convergence rate and key parameters. IBGWO utilizes binary values for wolf positioning and combines distance metrics, while SVM classifies data effectively, achieving high accuracy and low error rates.

Almaslukh et al. [23] proposed a deep learning-based intrusion detection method for WSNs with four stages. First, data is preprocessed and grouped using K-Means and DBSCAN clustering. Then, features like mean and variance are extracted and converted into embedded vectors with Word2Vec for CNN processing. This approach achieves high accuracy and rapid detection, enhancing intrusion detection effectiveness in WSNs.

Jingjing et al. [24] introduced the MC-GRU neural network, which uses dual-layer gradients for forward and backward propagation, with preprocessing to enhance detection accuracy. Particle swarm optimization fine-tunes network weights, improving performance. This method offers high accuracy and adaptability but depends on large datasets and is sensitive to noise, presenting some limitations.

Sood et al. [25] proposed a method using Conditional Generative Adversarial Networks (CGAN) to enhance intrusion detection by generating synthetic data. Steps include data preprocessing via wavelet transform, feature extraction with CNNs, and conditional data generation. This approach offers high accuracy and privacy benefits but is complex and demands extensive training data.

Zhang [26] proposed a credit scoring system for detecting intrusions in WSNs. It calculates credit scores based on node behavior and uses a Bayesian estimation algorithm to identify suspicious nodes. While the method provides good accuracy, it is computationally complex and may be time-consuming.

Gautami et al. [27] introduced the AQN3 intrusion detection system for WSNs. Data is collected using the DFFF algorithm, followed by preprocessing, matrix reduction with ELDA, and classification with QN3. This method accurately detects intrusions and various attacks but requires substantial processing resources.

Wu et al. [28] introduces Bot-DM, a novel dual-mode botnet detection method leveraging raw network traffic. Unlike traditional approaches that rely on deterministic features, Bot-DM analyzes both the payload and header information to detect abnormal behaviors. It employs a multi-layer Transformer encoder to extract implicit semantic relationships in payload data (modelled as text) and uses graphical representations of headers to capture spatial features. By maximizing the mutual information between these two modalities, Bot-DM enhances detection accuracy for both known and unknown botnets. Experimental results demonstrate superior performance, achieving up to 99.84% accuracy on public datasets, significantly outperforming existing methods.

Ali et al. [29] addresses the critical challenge of detecting botnets in IoT environments by proposing a hybrid deep learning model that combines LSTM Autoencoders and Multilayer Perceptrons (MLPs). This fusion enables effective analysis of sequential data and intricate pattern recognition, enhancing the detection of botnet activities and zero-day attacks. Evaluated on two large IoT traffic datasets, N-BAIoT2018 and UNSW-NB15, the model achieved impressive accuracy rates of 99.77% and 99.67%, respectively, outperforming existing detection systems. Key contributions include its robustness against zero-day attacks and its validation through extensive real-world testing. Despite its promising results, challenges such as server security, scalability, and training complexity must be addressed to ensure its broader applicability and reliability in industrial-scale IoT security.

A summary and comparison of the mentioned methods is present in Table 1.

**Table 1.** Summary and comparison of related works

| Ref | Year | Purpose | Key features | Advantages | Disadvantages |
|---|---|---|---|---|---|
| [15] | 2017 | Intrusion detection based on trust and energy | Using Kalman filter for trust prediction and double ACK algorithm for authentication | Improve packet delivery rate and reduce delay, Reducing vulnerabilities | High control overhead and more energy consumption than MS-MAC |
| [16] | 2018 | Intrusion Detection with Structural Adaptation Capability | Unsupervised Learning and Mean Shift Clustering Algorithm for Intrusion Detection, Weighted Support Vector Machine | High Detection Rate with Low False Alarm Rate, Adaptability to Network Structural Changes | High Processing Overhead |
| [17] | 2018 | Intrusion Detection with Intelligent Optimization and Fuzzy Clustering | Feature Optimization using Enhanced Social Spider Algorithm and Distributed Fuzzy Clustering, Parallel Processing | Improved Performance in Large-Scale Networks through Parallel Processing and Fuzzy Clustering | Not yet deployed in real-world networks; Insufficient emphasis on success rate with small sample sizes |
| [18] | 2019 | Trust-based intrusion detection at the protocol layer | Identification of different intrusions in different layers, Reducing memory and energy overhead | High accuracy, reducing energy consumption and complexity of messages | Lack of evaluation in real environments with outdoor modules |

| Ref | Year | Purpose | Key features | Advantages | Disadvantages |
|---|---|---|---|---|---|
| [19] | 2019 | Diagnosis of short-term and long-term abnormalities | Integration of short-term and long-term approaches for accurate detection of anomalies | high accuracy and timing of anomalies; Reduction of false positives | Need to improve process automation and better set thresholds |
| [20] | 2020 | Intrusion detection for mobile networks | Using dedicated sniffers and linear regression; Cross-layer architecture | High accuracy by reducing misleading data and eliminating false positives | Need to improve threshold settings and performance of sniffers in complex environments |
| [21] | 2020 | Hierarchical intrusion detection in wireless networks | combination of REP Tree, JRip and Forest PA classifiers; Reduce the impact of data imbalance | Minimizing the impact of data imbalance, increasing accuracy and preventing calculation errors | Requires high computational resources for multi-class combination |
| [22] | 2021 | Improving the accuracy of intrusion detection in wireless sensor network | Using IBGWO for higher convergence rate and SVM for intrusion detection | High accuracy, low error rate, reliable performance and high efficiency | The need to determine parameters, limitations in application, more processing time |
| [23] | 2021 | Intrusion detection in wireless sensor network with deep learning | Data preprocessing, feature extraction, entity embedding, neural network modeling | High accuracy and speed, cost reduction, reusability | The need for big data, powerful hardware, complex engineering |
| [24] | 2022 | Intrusion detection in wireless sensor networks with high accuracy and efficiency | Using recurrent neural network with direct and inverse gradient layers, data preprocessing with PCA algorithm and dimensionality reduction | High accuracy due to the use of MC-GRU, high processing speed and improved efficiency, adaptability to a large number of sensors | Requirement of large training data, sensitivity to noise and heterogeneous inputs |
| [25] | 2022 | Improving intrusion detection accuracy in wireless sensor networks by generating conditional data | Using Conditional Generative Neural Network (CGAN), generating conditional data to increase accuracy, preprocessing data with wavelet transform | High accuracy using CGAN conditional data, high adaptability and model adjustment based on network conditions, privacy protection using fake data | Need more data for training, high computational complexity and cost, sensitivity to noise and erroneous data |
| [26] | 2023 | Detection of malicious nodes in the network with credibility score | Calculating the credit score of nodes, using the trust improvement system | Good accuracy, reliable for nodes | Being time-consuming, requiring complex calculations |
| [27] | 2023 | Intrusion detection with AQN3 system in wireless sensor network | Cluster leader selection with DFFF, classification of attack types, matrix reduction with ELDA | High accuracy in detecting types of attacks, high adaptability | The need for high processing resources, high computational cost |
| [28] | 2024 | To propose a novel botnet detection method, Bot-DM, that leverages dual-modal analysis of raw network traffic. | Combines payload (text-based) and header (graphical) data. | High detection accuracy No dependence on predefined traffic field meanings. Adaptable to complex and emerging botnet scenarios. | Requires significant computational resources. Assumes availability of large and diverse datasets for training. |
| [29] | 2024 | To develop a robust hybrid deep learning model for detecting botnets and zero-day attacks in IoT networks. | Combines LSTM Autoencoders for sequential data analysis with MLP for pattern recognition. | High accuracy rates Robust against zero-day attacks by identifying unusual patterns. Effective for real-world IoT scenarios. | Reliance on a central server introduces potential security risks. Scalability challenges for industrial applications. High training complexity and dependence on large datasets. |

## 3. Problem Statement

As the Internet of Things continues to expand across various sectors—such as healthcare, education, agriculture, and smart cities—IoT networks face significant security challenges. The rapid proliferation of IoT devices has created vast, interconnected networks with heterogeneous devices, complex data flows, and constrained resources. This growth has also made IoT networks increasingly vulnerable to various forms of cyberattacks, including unauthorized access, data breaches, and denial-of-service (DoS) attacks, which can disrupt services, compromise sensitive information, and lead to significant operational risks.

While traditional security measures like encryption, authentication, and access control form the initial defense for IoT networks, these preventative measures alone are insufficient to counter increasingly sophisticated and adaptive cyber threats. Once these defenses are bypassed, it becomes crucial to detect intrusions in real time to protect the network. Intrusion Detection Systems (IDS) have thus become an essential component in securing IoT networks, providing a second layer of defense by identifying unusual activities and alerting network administrators to potential threats. However, existing IDS solutions are often challenged by the unique requirements of IoT environments. These systems must achieve high detection accuracy while operating within the resource limitations typical of IoT devices, such as low processing power, limited memory, and constrained energy. Additionally, achieving high accuracy without excessive false positives remains a critical goal, as unnecessary alerts can drain resources and reduce system reliability.

Given these constraints, there is a need for an IDS capable of effectively identifying intrusions in IoT networks with a balance between accuracy and resource efficiency. The proposed research seeks to address this gap by developing a two-phase intrusion detection model specifically tailored for IoT. In the first phase, the Grasshopper Optimization Algorithm (GOA) is used for feature selection, reducing the dataset to its most relevant features to minimize processing load and improve model efficiency. In the second phase, a Support Vector Machine classifier is employed to accurately categorize network traffic into normal or malicious activity.

This dual-phase approach is expected to enhance detection accuracy and reduce false positives while maintaining low computational overhead, making it highly suited to IoT environments. By achieving these objectives, the proposed model aims to provide a scalable and efficient solution to strengthen IoT network security, ensuring that connected devices can operate reliably and securely in increasingly complex and dynamic settings.

## 4. Proposed Method

Security management in Internet of Things networks typically involves two primary approaches: prevention-based and detection-based strategies. In the proposed method, detection-based management is employed. In an IoT security framework, if prevention techniques—such as encryption, authentication, and validation, which serve as the first line of defense—are bypassed, the focus then shifts to the second line of defense: intrusion detection. Intrusion detection is critical for the alert mechanism within IoT systems, enabling effective security monitoring across interconnected devices. An intrusion detection system (IDS) in IoT can be either a software or hardware tool that actively monitors network traffic to identify both internal and external threats. IDSs analyze system and device activities, recognize known attack patterns, and detect abnormal behaviors, all essential for managing the security of diverse IoT environments.

In the proposed method, the Grasshopper Optimization Algorithm (GOA) is used for feature reduction, which helps select the most relevant features for accurate and efficient detection. This is especially beneficial in IoT settings, where optimizing resource usage is crucial. Following this, a Support Vector Machine (SVM) classifier is applied to categorize data samples, effectively distinguishing between normal and potentially malicious activities. This two-step approach enhances detection accuracy and processing efficiency, meeting the unique requirements of IoT security management. An overview of the general steps in the proposed method is depicted in Figure 1.

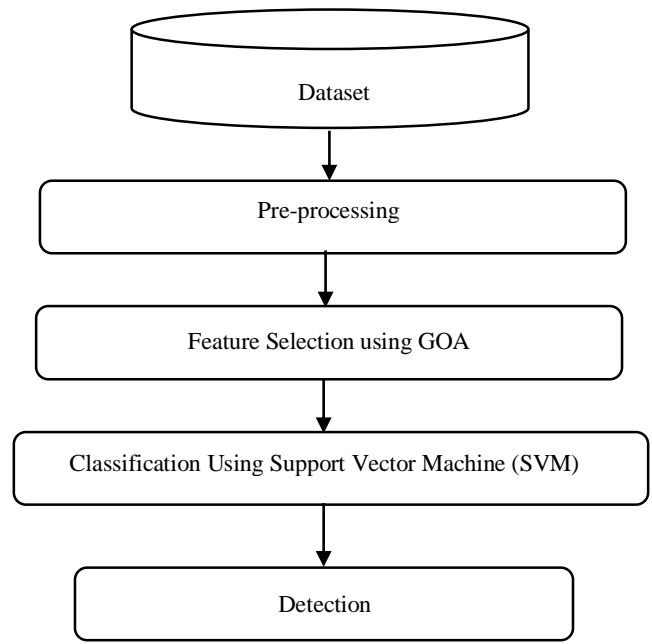

**Figure 1**: Steps of the proposed method

Alt Text Figure 1[54 words]: The image shows a flowchart illustrating the steps of a proposed method. It starts with a "Dataset" at the top, followed by "Pre-processing," then "Feature Selection using GOA," next is "Classification Using Support Vector Machine (SVM)," and finally "Detection" at the bottom. The figure is labeled as "Figure 1: Steps of the proposed method."

### 4-1. Pre-processing

In an Internet of Things (IoT) network, datasets often contain noise, redundancy, and various data types, including numerical and string values. Extracting meaningful insights and building models from these datasets presents several challenges. For example, string values cannot be directly utilized in algorithms that require numerical inputs, necessitating conversion to a numerical format. The process for converting string data to numerical values is outlined in Algorithm 1.

| Algorithm. 1: Converting String Data to Numeric Values |
|---|
| **Input:** Dataset $D$ <br> **Output:** processed Dataset $D'$ <br> (1) **for** $n = 1$ to $N$ do <br> (2)     **for** $m = 1$ to $N$ <br> (3)         $(S, c) \leftarrow \text{find}(D)$ |

```
(4)     end for
(5) end for
(6) for n = 1 to c do
(7)     D'(S) ← D(m)
(8) end for
(9) return D'
```

To convert string data to numbers, the dataset DDD is examined, and all string entries SSS within it are identified. The find() function is used to locate the corresponding columns, denoted by ccc, as shown in lines 1 to 3 of Algorithm 1. Subsequently, the replace() function is applied to substitute each string value SSS with a random numerical value mmm, as indicated in lines 7 to 9 of the algorithm. Finally, the processed dataset D′D'D′ is generated, ready to serve as the input for the next step. This pre-processed dataset, now consisting of numerical values, is suitable for use in intrusion detection algorithms.

In this paper, the NSL-KDD dataset is used as the input, comprising 125,973 records, 41 features, and 5 classes. This dataset includes a normal class and four attack classes, specifically: DoS, U2R, R2L, and Prob.

### 4-2. Feature Selection

The Grasshopper Optimization Algorithm (GOA) is used to perform feature selection, identifying the most critical features within the dataset. The goal is to categorize features into two groups: significant features that influence attack detection and non-significant features that have minimal impact on this process.

The feature selection process using the Grasshopper Optimization Algorithm consists of the following steps:

- **Step 1:** Initialize parameters.
- **Step 2:** Evaluate the population of grasshoppers and select or update TTT.
- **Step 3:** Check the stop condition.
- **Step 4:** Update the location of the grasshoppers.
- **Step 5:** Update parameter c.
- **Repeat**: Return to Step 2.

The flowchart of the proposed feature selection method using the Grasshopper Optimization Algorithm is presented in Figure 2.

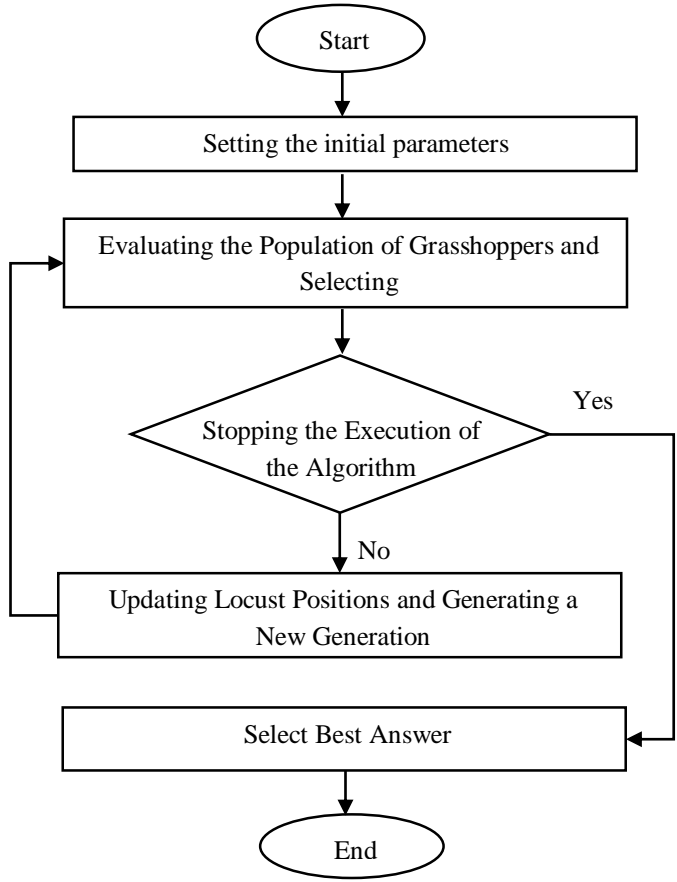

**Figure 2**: Feature selection flowchart in the proposed method

Alt Text Figure2[75 words]: The image is a flowchart titled "Feature selection flowchart in the proposed method." It starts with "Start" and proceeds to "Setting the initial parameters." Then it moves to "Evaluating the Population of Grasshoppers and Selecting." A decision diamond follows, labeled "Stopping the Execution of the Algorithm," with "Yes" leading to "Select Best Answer" and then "End." If "No," it goes to "Updating Locust Positions and Generating a New Generation," looping back to the evaluation step.

**Step 1: Setting the Initial Parameters**

To address the feature selection problem for intrusion detection in an Internet of Things (IoT) environment, a solution based on the Grasshopper Optimization Algorithm (GOA) is applied, starting with the design of the problem's search space. In this space, grasshoppers (agents) explore potential solutions.

For the NSL-KDD dataset used in this research, the initial population of GOA is represented as a set of 41-element binary arrays, where each element corresponds to a specific attribute in the dataset. A value of "1" indicates that the feature is significant for detecting intrusions in the IoT network, while a value of "0" indicates that the feature is unimportant or ineffective in identifying attack patterns. At this stage, initial solutions are encoded in binary form, and the initial population is randomly generated.

To generate the initial population, each element in the 41-element array is assigned a random value of either "1" or "0." For instance, in a 41-element array, 20 elements might be randomly set to "1" (indicating selected features), and the remaining elements are set to "0" (non-selected features). This initial population represents various solutions, each positioned as a point within the problem's search space.

As illustrated in Figure 3, the initial population in the Grasshopper Optimization Algorithm is generated using binary encoding. Each element of this population constitutes a potential solution to the optimal feature selection problem.

|  | Feature1 | Feature2 | Feature3 | Feature4 | . | . | . | Feature39 | Feature40 | Feature41 |
|---|---|---|---|---|---|---|---|---|---|---|
| Solution1 | 1 | 0 | 0 | 1 | 0 | 0 | 0 | 1 | 1 | 1 |
| Solution2 | 0 | 0 | 1 | 1 | 0 | 0 | 0 | 0 | 1 | . |
| Solution3 | 1 | 0 | 1 | 0 | 0 | 0 | 0 | 0 | 0 | 0 |
| | . | . | . | . | . | . | . | . | . | . |
| | . | . | . | . | . | . | . | . | . | . |
| | . | . | . | . | . | . | . | . | . | . |
| Solution n-1 | 1 | 1 | 0 | 0 | 0 | 0 | 0 | 1 | 1 | 0 |
| Solution n | 1 | 1 | 0 | 0 | 0 | 0 | 0 | 1 | 1 | 0 |

**Figure 3**: An example of the initial population

Alt Text Figure3[74 words]: The image is a table labeled "Figure 3: An example of the initial population." It shows a binary representation of different solutions (Solution1 to Solution n) across various features (Feature1 to Feature41). Each solution consists of a series of 0s and 1s, indicating the presence (1) or absence (0) of a feature in that solution. The table illustrates how different combinations of features are selected to form an initial population for a computational method.

In the Grasshopper Optimization Algorithm, each solution array represents a subset of selected features. For instance, in Solution 1, the attributes in positions 1, 4, ..., 39, 40, and 41 are considered significant (indicated by a value of "1") and are therefore used in this solution.

**Step 2: Evaluating the Population of Grasshoppers and Selecting/Updating T**

During the evaluation phase, the initial population generated in the previous step is assessed. Each initial solution, represented as a point within the search space, is evaluated using a performance index to identify the most promising solutions for intrusion detection in an Internet of Things.

The performance index in this phase can include the following metrics:
- Accuracy: The ratio of correct predictions to the total number of predictions.
- Sensitivity: The ratio of true positive predictions to the total number of positive samples.
- Selected Features: The count of features used in the model.
- Speed: The time taken to complete the prediction.

After evaluating the initial population, the top-performing solutions are selected to continue generating the next generation of solutions. The primary goal is to determine the optimal or near-optimal set of features that maximizes detection accuracy while minimizing errors.

In this step, each solution array is assessed to determine which features have the greatest impact on intrusion detection accuracy and error reduction. The solution with the highest fitness function value is selected as the optimal solution.

The fitness function is a key component of the Grasshopper Optimization Algorithm and determines whether a solution is viable. It should minimize the error rate, maximize the true positive rate, and select an effective subset of features. In the proposed method, the fitness function evaluates each feature subset using three key parameters: the true positive rate, the error rate, and the number of selected features. The fitness function is calculated as shown in Equation (1).

$$Fitness = R_{tp} + (1 - R_E) + (1 - \frac{N_F}{41}) \tag{1}$$

In Equation (1), $R_{tp}$ represents the true positive rate, calculated according to Equation (2), $R_E$ denotes the error rate, calculated as per Equation (3), and $N_F$ is the number of selected features. Based on Equation (1), a higher correct positive rate, a lower error rate, and a smaller feature set all contribute to maximizing the fitness function.

The true positive rate $R_{tp}$ is calculated by dividing the number of true positive cases (TP) by the sum of true positives (TP) and false negatives (FN). This calculation is shown in Equation (2):

$$R_{tp} = \frac{TP}{TP + FN} \tag{2}$$

Where TP is the count of true positive cases, and FN is the count of false negative cases.

The error rate $R_E$ represents the proportion of cases that the classifier incorrectly predicts out of all cases. It is calculated using Equation (3):

$$R_E = \frac{FP + FN}{TP + FN + TN + FP} \tag{3}$$

Where TP is the number of true positives, FP is the number of false positives, TN is the number of true negatives, and FN is the number of false negatives.

These metrics collectively assess the performance of each solution, aiming to maximize detection accuracy and minimize error rates while using the smallest feature subset possible.

**Step 3: Stopping the Execution of the Algorithm**
The criteria for stopping the grasshopper optimization algorithm in the proposed method are as follows:
- The algorithm reaches a maximum of 40 iterations.
- The difference in the maximum fitness between the current and previous generations is less than 0.001.

If either condition is met, the algorithm stops. Upon stopping, the best solution from the population is selected for determining the optimal features for intrusion detection. If neither condition is met, the algorithm proceeds to the next step.

**Step 4: Updating Locust Positions and Generating a New Generation**
The new position of each locust is calculated using Equation (4):

$$X_i = \sum_{j=1, j \neq i}^{N} s(|x(i) - x(j)|) \frac{|x(i) - x(j)|}{d_{ij}} - g\widehat{e_g} + u\widehat{e_w} \tag{4}$$

Where u is the thrust constant, and $\widehat{e_w}$ is the unit vector in the wind direction. g represents gravitational force, and $\widehat{e_g}$ is the unit vector toward the Earth's center. $d_{ij}$ is the distance between locust i and locust j. $X_i$ and $X_j$ are the positions of locusts i and j, respectively, and N is the total number of locusts.

Each locust's next position is determined using Equation (5):

$$X_i^d(t+1) = c \left\{ \sum_{j=1, j \neq i}^{N} c \frac{ub_d - lb_d}{2} s(|x(i) - x(j)|) \frac{x(i) - x(j)}{d_{ij}(t)} \right\} + \widehat{T_d(t)} \tag{5}$$

Where $ub_d$ and $lb_d$ are the upper and lower limits in dimension d. $\widehat{T_d(t)}$ represents the best solution up to iteration t in dimension d. c is the reduction factor.

Two additional operations, SWAP and Reversion, are employed to improve solution diversity:
- **SWAP:** Two features are randomly selected, and their positions are exchanged.
- **Reversion:** Two features are selected, and all features between them are inverted.

This process generates a new population of solutions, iterating towards the optimal feature selection for the intrusion detection model.

**Step 5: Updating Parameter c**
In this step, parameter c is adjusted to control the influence of the locusts' movement as the algorithm progresses. This adjustment is achieved using Equation (6):

$$c = c_{max} - t \frac{c_{max} - c_{min}}{tmax} \tag{6}$$

Where $c_{max}$ is the initial maximum value of c, typically close to 1. $c_{min}$ is the minimum value of c, usually near zero. t is the current iteration number. $tmax$ is the maximum number of iterations. As the algorithm iterates, c gradually decreases from $c_{max}$ to $c_{min}$, controlling how much influence each locust's position has on others. This approach helps balance the exploration and exploitation phases of the optimization, allowing for broad searching initially and more precise adjustments in later iterations.

After updating c, the algorithm returns to the evaluation stage (Step 2), where the solutions are re-evaluated, and the process continues until the stopping criteria are met.

### 4-3. Classification Using Support Vector Machine (SVM)

Once the optimal features have been selected, the data is classified using a Support Vector Machine (SVM). SVM is a powerful machine learning algorithm that can handle both linear and nonlinear data classification. By transforming the data into a higher-dimensional space using a nonlinear mapping, SVM finds a hypersurface that effectively separates examples of each class. For two-class classification, SVM identifies a hyperplane that maximizes the margin, or the distance, between data points of each class. In cases where there are multiple classes (e.g., N classes), the SVM implements a one-vs-all strategy:

- Each class is trained against all other classes, where the data of the target class is labeled as +1, and all other classes are labeled as -1.
- This approach trains **N support vector machines**, each corresponding to one class.
- During testing, each sample is evaluated by all N SVMs, and the sample is assigned to the class with the highest output.

Steps for Classification with SVM are as follows:

**Step 1: Define Data Points**

The dataset is represented as $\{(x_1.c_1).(x_2.c_2).\ldots.(x_n.c_n)\}$, where each $x_i$ is a p-dimensional vector representing the selected features, and each $c_i$ is the class label for $x_i$

**Step 2: Define Parallel Hyperplanes**

SVM finds two parallel hyperplanes in the feature space that maximize the margin between classes. This approach selects a separator with the maximum distance from the closest data points of each class. These closest data points are known as support vectors**.** The equations of these two parallel hyperplanes are defined as Equation (7) and Equation (8):

$$w.x - b = 1 \qquad (7)$$
$$w.x - b = -1 \qquad (8)$$

Where w represents the weight vector perpendicular to the hyperplane. b is a bias term**.**

**Step 3: Determine the Optimal Separator**

The separator lies equidistantly between these two planes. The distance between them is maximized to enhance the classifier's generalization ability.

By classifying data points based on their position relative to these hyperplanes, the Support Vector Machine (SVM) effectively distinguishes between classes, achieving a high level of accuracy in differentiating between various attack types and normal traffic within an Internet of Things (IoT) network. Given the heterogeneous and resource-constrained nature of IoT devices, this precision is particularly valuable, as it enables reliable identification of potential security threats with minimal computational overhead. In IoT networks, where data flows continuously across numerous connected devices, the SVM classifier can play a pivotal role in maintaining security by quickly recognizing abnormal patterns indicative of cyberattacks. This capability not only supports timely threat detection but also enhances the overall resilience of the IoT ecosystem by adapting to evolving attack vectors and effectively managing the high volume and variety of data typical in IoT environments.

In Figure (4), two groups of data are shown, each belonging to a different class. The SVM algorithm's task is to find the hyperplanes (lines in this case) that separate these groups with the largest margin. The margin is the shortest distance between the hyperplanes and the closest data points from both classes, called support vectors.

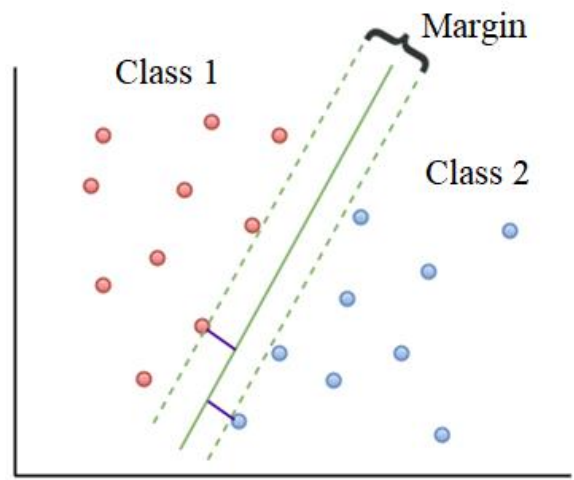

**Figure 4:** linear classification hyperplanes in 2D space

Alt Text Figure4[53 words]: The image depicts a two-dimensional plot of a linear classification model. The red and blue points represent class 1 and class 2, respectively. A green line serves as the decision boundary between the two classes. Two light green parallel lines indicate the margins, representing the maximum distance between them and the decision boundary.

While there are many possible hyperplanes that could separate the data points, the SVM algorithm selects the one that maximizes the margin. If the training data is linearly separable, two boundary hyperplanes can be chosen such that no data points lie between them, and the margin is maximized. The margin between the two hyperplanes is given by the Equation (9):

$$\text{Margin} = \frac{2}{|w|} \tag{9}$$

Where w is the weight vector, which is perpendicular to the hyperplane. The goal is to minimize |w| to maximize the margin.

To ensure that no data points fall within the margin, mathematical constraints are added to the model. For each data point iii, the following constraints must hold:

For data points of the first class (positive class):

$$w.x_i - b \geq 1 \tag{10}$$

For data points of the second class (negative class):

$$w.x_i - b \leq -1 \tag{11}$$

Where w is the weight vector. $x_i$ is the feature vector of data point i. b is the bias term.

These constraints ensure that each data point is correctly classified and does not lie within the boundary region between the two hyperplanes.

After applying the feature selection process, the SVM weights are assigned based on the optimal subset of features. These weights are used to train the SVM, and the classification accuracy is evaluated. The chromosome (or solution) with the highest classification accuracy is selected as the optimal solution.

Once the optimal SVM parameters (weights and support vectors) are identified, they are decoded to form the final intrusion detection model.

The objective of proposed method is to propose a novel intrusion detection model capable of effectively identifying malicious nodes within Internet of Things networks. Given the unique characteristics of IoT, including diverse device types, limited resources, and extensive data exchange, this model leverages Support Vector Machine to classify network traffic into two primary categories:

- Normal traffic
- Abnormal traffic (indicating a potential attack)

The model employs a feature selection process to identify the most relevant attributes, enhancing both detection accuracy and computational efficiency—key factors in IoT environments. Following this optimized feature selection, SVM is used to classify data, enabling precise and timely identification of intrusion attempts. This approach aims to bolster IoT security by reliably detecting abnormal activities while conserving device resources, ultimately improving the overall resilience and safety of IoT networks against a range of cyber threats.

## 5. Simulation and Evaluation

The proposed method was simulated using MATLAB 2017, which provides an efficient platform for implementing and evaluating machine learning algorithms, including feature selection and classification. The simulation parameters have been provided in Table 2.

**Table 2:** simulation parameters

| Parameter | Value |
|---|---|
| Simulation time | 1000 s |
| package size | 4 kb |
| Network size | 100*100 m$^2$ |
| The number of nodes | 100 |
| Primary energy | 0.5 J |

## 5-1. Dataset

For the simulation, the NSL-KDD dataset is used, which is a corrected and enhanced version of the KDD99 dataset. Both datasets were originally created by MIT Labs and DARPA in 1998, with the goal of simulating normal and abnormal communication behaviour in military network environments. However, the KDD99 dataset contained inherent flaws, such as class imbalance and redundancy, which reduced the accuracy of intrusion detection systems. These flaws were addressed in the NSL-KDD dataset.

- NSL-KDD is a more balanced version of the KDD99 dataset and is considered a more realistic benchmark for evaluating intrusion detection systems.
- NSL-KDD is designed to provide a challenging evaluation environment because it consists of more difficult instances identified using clustering methods. This makes it a suitable dataset for testing the effectiveness of intrusion detection systems, ensuring that they are not biased toward easily detectable attacks.
- NSL-KDD contains 41 features for each sample, representing various network and host-level attributes. Each record also includes a label indicating whether the connection is normal or represents one of four types of attacks.

The dataset categorizes attacks into four types:

- Denial of Service (DoS): Description: These attacks overwhelm a victim's resources, preventing legitimate access or services. Related Attributes can include "Source bytes" and "Percentage of packets with errors."
- Probe: Description: These attacks are reconnaissance in nature, where the attacker tries to gather information about the target system, such as through port scanning. Related Attributes can include "Connection duration" and "Resource bytes."
- User-to-Root (U2R): Description: The attacker gains access to a user account and attempts to escalate privileges to the root level using vulnerabilities such as buffer overflow. Related Attributes can include "Number of files created" and "Number of requests considered."
- Root-to-Local (R2L): Description: These attacks occur when an attacker gains access to a remote machine and attempts to gain local access by exploiting vulnerabilities like password guessing. Related Attributes can include "Connection time" and "Service request" at the network level, and "Number of failed login attempts" at the host level.

Each connection in the dataset is labeled as normal or abnormal (depending on the attack type). The dataset contains both continuous (e.g., byte counts) and discrete features (e.g., attack labels), making it a complex challenge for intrusion detection models.

### 5-2. Evaluation parameters

In this research, the distributed penetration detection framework is evaluated based on detection performance and false positive rates. These metrics provide insight into how effectively the system identifies malicious activities and how often it misclassifies normal activities as attacks.

The performance of the intrusion detection system (IDS) is evaluated using several classification metrics, which are as follows:

- **True Positive (TP)**: Definition: The true class of the test sample is positive, and the classifier correctly classifies it as positive.
- **False Negative (FN)**: Definition: The true class of the test sample is positive, but the classifier incorrectly predicts the class as negative.
- **False Positive (FP)**: Definition: The true class of the test sample is negative (normal), but the classifier incorrectly classifies it as positive (attack).
- **True Negative (TN)**: Definition: The true class of the test sample is negative (normal), and the classifier correctly identifies it as negative.

Based on the above definitions, the following evaluation parameters are calculated:

True Positive Rate (TPR), also known as Sensitivity or Recall: This metric indicates how effective the system is at identifying attacks. A higher TPR means the IDS is good at detecting malicious behaviours.

$$TPR = \frac{TP}{TP + FN} \tag{12}$$

False Positive Rate (FPR): This metric measures the number of normal activities misclassified as attacks. A lower FPR is desirable to avoid unnecessary alarms.

$$FPR = \frac{FP}{FP + TN} \tag{13}$$

True Negative Rate (TNR): A high TNR means that the system is also good at recognizing normal behavior, not just attacks.

$$TNR = \frac{TN}{TP+FN} \tag{14}$$

False Negative Rate (FNR): A low FNR is desirable as it indicates that fewer attacks are missed and incorrectly classified as normal.

$$FNR = \frac{FN}{FN + TP} \tag{15}$$

Accuracy: This metric gives a general idea of the system's overall performance. However, accuracy alone may not be sufficient to evaluate performance in imbalanced datasets where attacks are much rarer than normal activities.

$$Accuracy = \frac{(TP+TN)}{(TP+FN+FP+TN)} \tag{16}$$

The goal of the proposed method is to maximize the True Positive Rate (TPR) and True Negative Rate (TNR) while minimizing the False Positive Rate (FPR) and False Negative Rate (FNR). The overall Accuracy will provide a general understanding of the IDS's effectiveness in detecting both attacks and normal connections.

These evaluation parameters provide a comprehensive view of the proposed method's performance and its ability to correctly classify different types of network behaviors. The next step in the evaluation process involves running the model on the NSL-KDD dataset and calculating these metrics to determine how well the system performs in detecting and preventing network intrusions.

## 5-3. Evaluation Results

The results of the proposed method are presented according to the evaluation parameters discussed earlier, such as True Positive Rate (TPR), False Positive Rate (FPR), True Negative Rate (TNR), False Negative Rate (FNR), and Accuracy. The performance was evaluated using ten-fold cross-validation, a common technique to assess the generalization of machine learning models.

Figure 5 shows the results of the penetration detection model for the proposed method, demonstrating how well the system classifies network traffic into normal and attack categories using the Grasshopper Optimization Algorithm (GOA) combined with the Support Vector Machine (SVM) classifier.

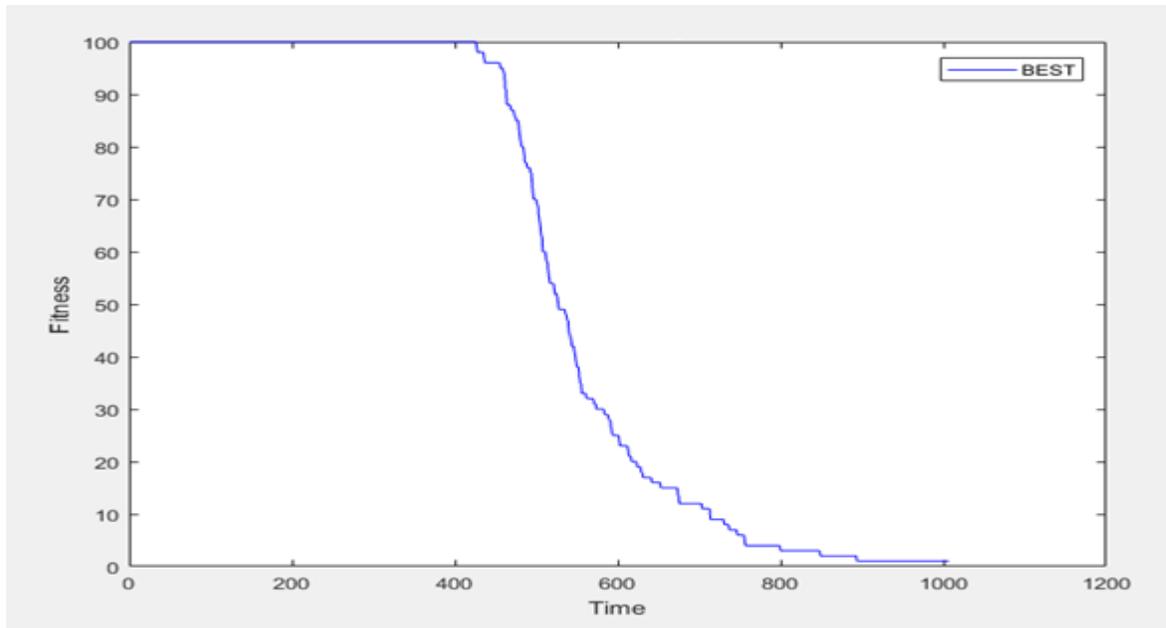

**Figure 5:** The output of the penetration detection model stage of the proposed method

Alt Text Figure5[60 words]: The image is a line graph titled "Figure 5: The output of the penetration detection model stage of the proposed method." The graph plots "Fitness" on the y-axis against "Time" on the x-axis. The line, labeled "BEST," starts high on the fitness scale and gradually decreases, showing a trend where fitness declines over time, eventually stabilizing at a lower value.

Figure 6 and Figure 7 show the detection rates for five key types of attacks (e.g., DoS, U2R, R2L, Probe, and Normal behaviour) in terms of TPR and FPR. The average results indicate that the proposed method performs well in detecting attacks while minimizing false positives.

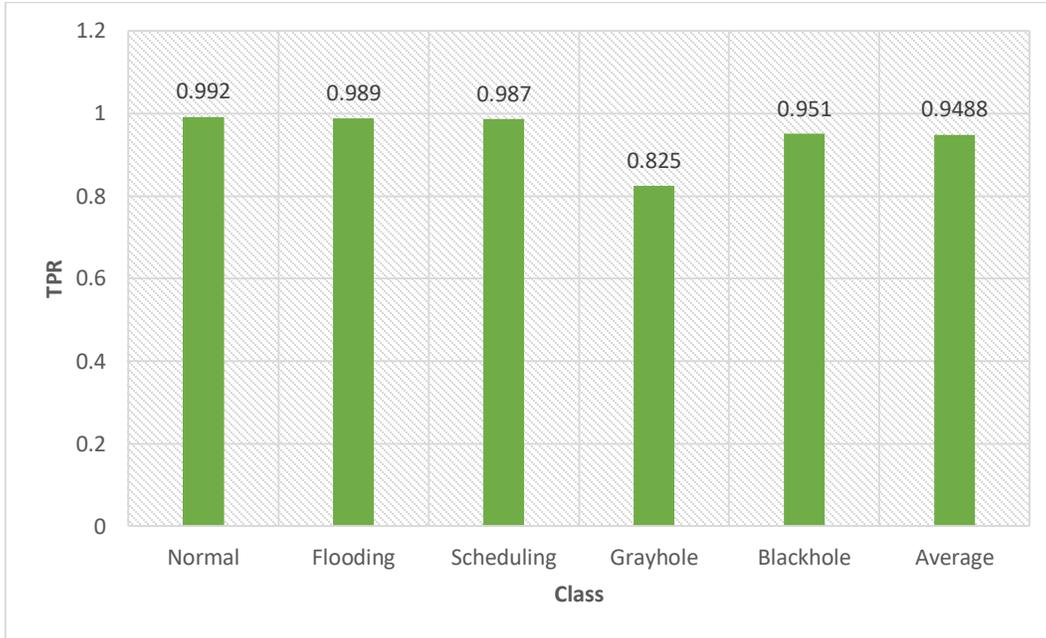

**Figure 6**: Comparison of the proposed method in terms of TPR parameter

Alt Text Figure 6[30 words]: Bar chart showing the TPR (True Positive Rate) of the proposed method across different classes. Values are: Normal (0.992), Flooding (0.989), Scheduling (0.987), Grayhole (0.825), Blackhole (0.951), and Average (0.9488).

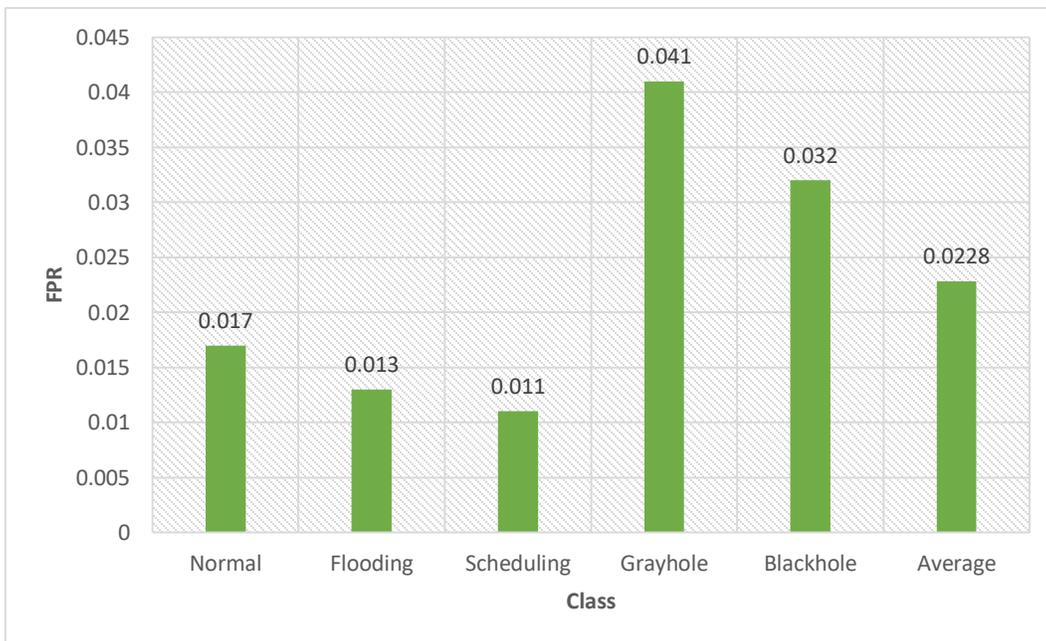

**Figure 7**: Comparison of the proposed method in terms of FPR parameter

Alt Text Figure 7[29 words]: Bar chart comparing the False Positive Rate (FPR) of the proposed method for six classes: Normal (0.017), Flooding (0.013), Scheduling (0.011), Grayhole (0.041), Blackhole (0.032), and the Average (0.0228).

Figure 8 and Figure 9 present the results for TNR and FNR, with the proposed method performing effectively in correctly classifying normal traffic and avoiding false negatives. This means fewer legitimate (normal) connections are misclassified as attacks, ensuring the model's reliability.

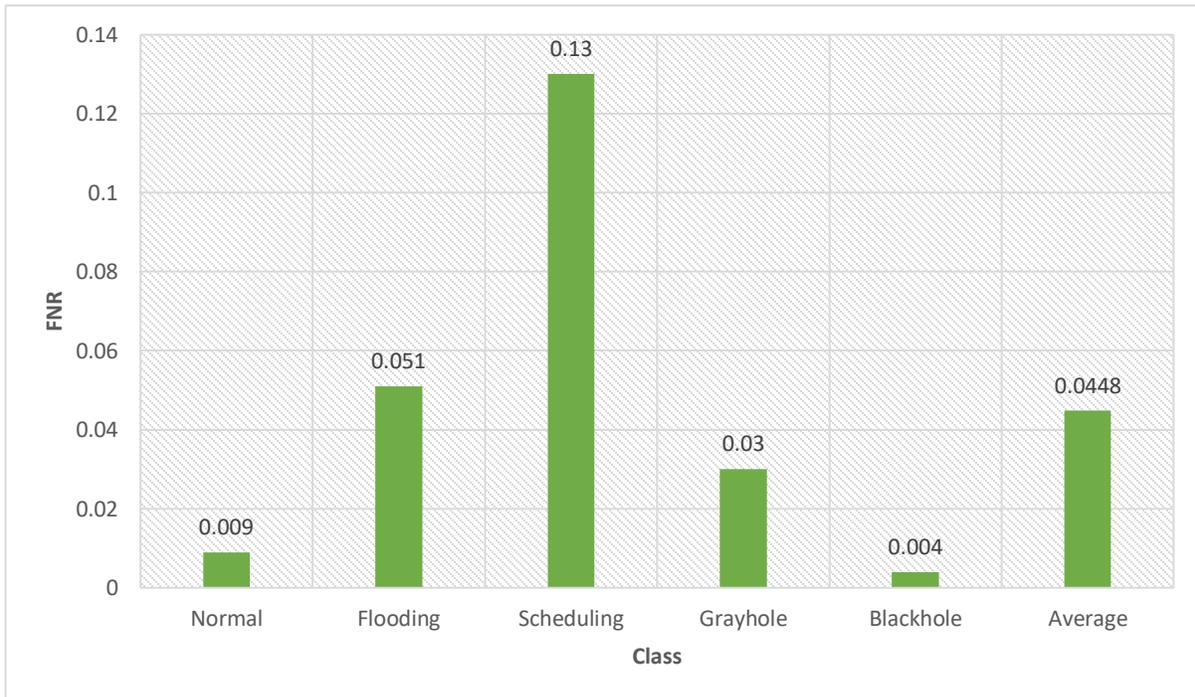

**Figure 8**: Comparison of the proposed method in terms of FNR parameter

Alt Text Figure 8[77 words]: The image shows a bar chart titled "Figure 8: Comparison of the proposed method in terms of FNR parameter." It compares the False Negative Rate (FNR) across different classes: Normal, Flooding, Scheduling, Grayhole, Blackhole, and Average. The Scheduling class has the highest FNR at 0.13, while the Blackhole class has the lowest at 0.004. Other classes have varying FNR values, with the Normal class at 0.009, Flooding at 0.051, Grayhole at 0.03, and an average of 0.0448.

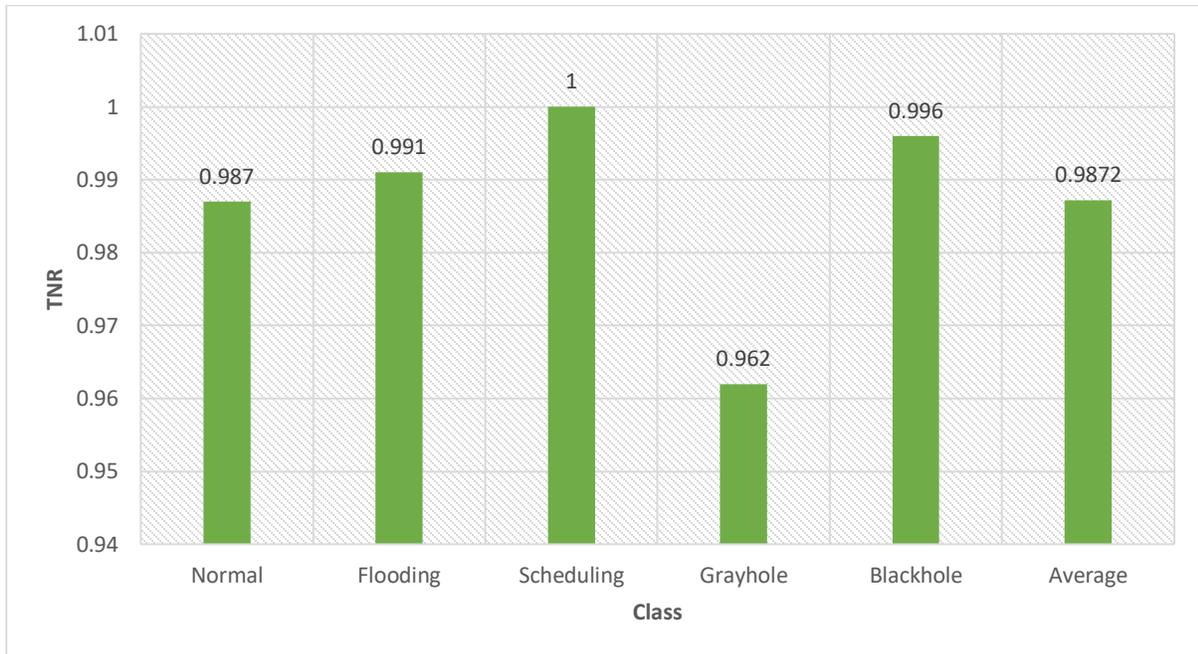

**Figure 9**: Comparison of the proposed method in terms of TNR parameter

Alt Text Figure9[29 words]: Bar chart comparing the True Negative Rate (TNR) of the proposed method for six classes: Normal (0.987), Flooding (0.991), Scheduling (1.000), Grayhole (0.962), Blackhole (0.996), and the Average (0.9872).

Figure 10 provides the classification accuracy for detecting the five important attacks and the overall average accuracy. The proposed method demonstrates high accuracy, showing its effectiveness in distinguishing between normal and attack traffic.

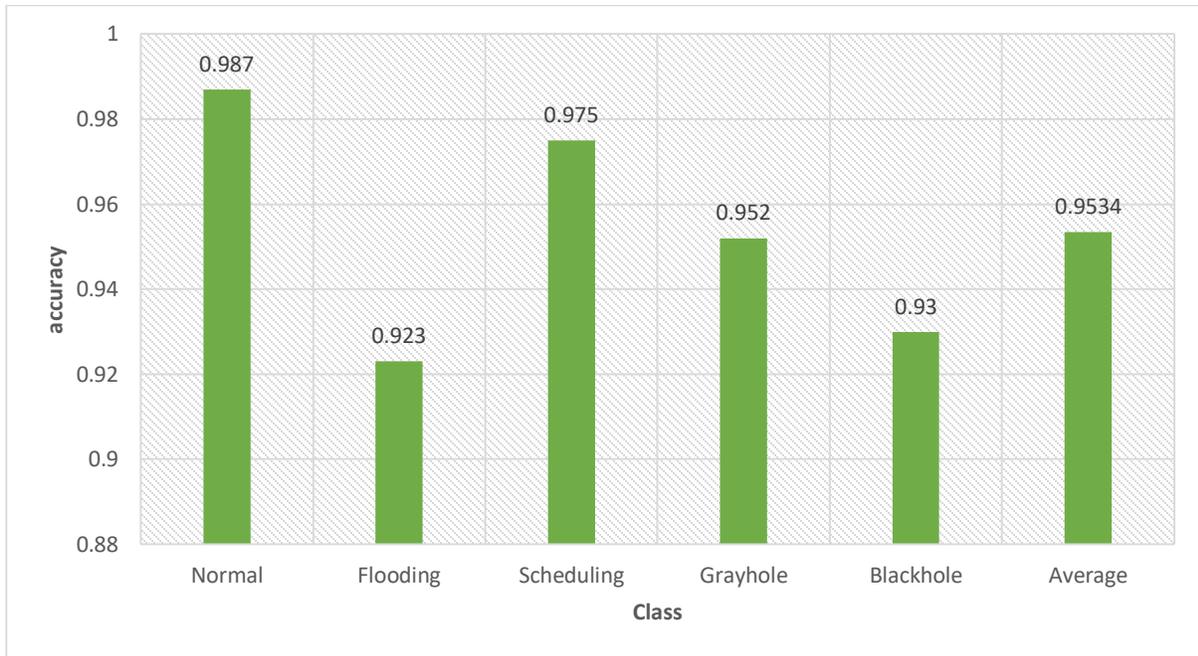

**Figure 10:** The result of the proposed method in terms of accuracy parameter

Alt Text Figure 10[27 words]: Bar chart showing the accuracy results of the proposed method for six classes: Normal (0.987), Flooding (0.923), Scheduling (0.975), Grayhole (0.952), Blackhole (0.930), and the Average (0.9534).

To further assess the effectiveness of the proposed method, it was compared with two other methods: Safaldin et al. [22] and Almaslukh et al. [23]. The results of these comparisons are presented below:

Figure 11 shows the comparison of the accuracy of the proposed method against the methods of [22] and [23]. The proposed method outperforms both methods, achieving a higher classification accuracy. This indicates that the integration of GOA for feature selection and SVM for classification yields better results compared to the other approaches.

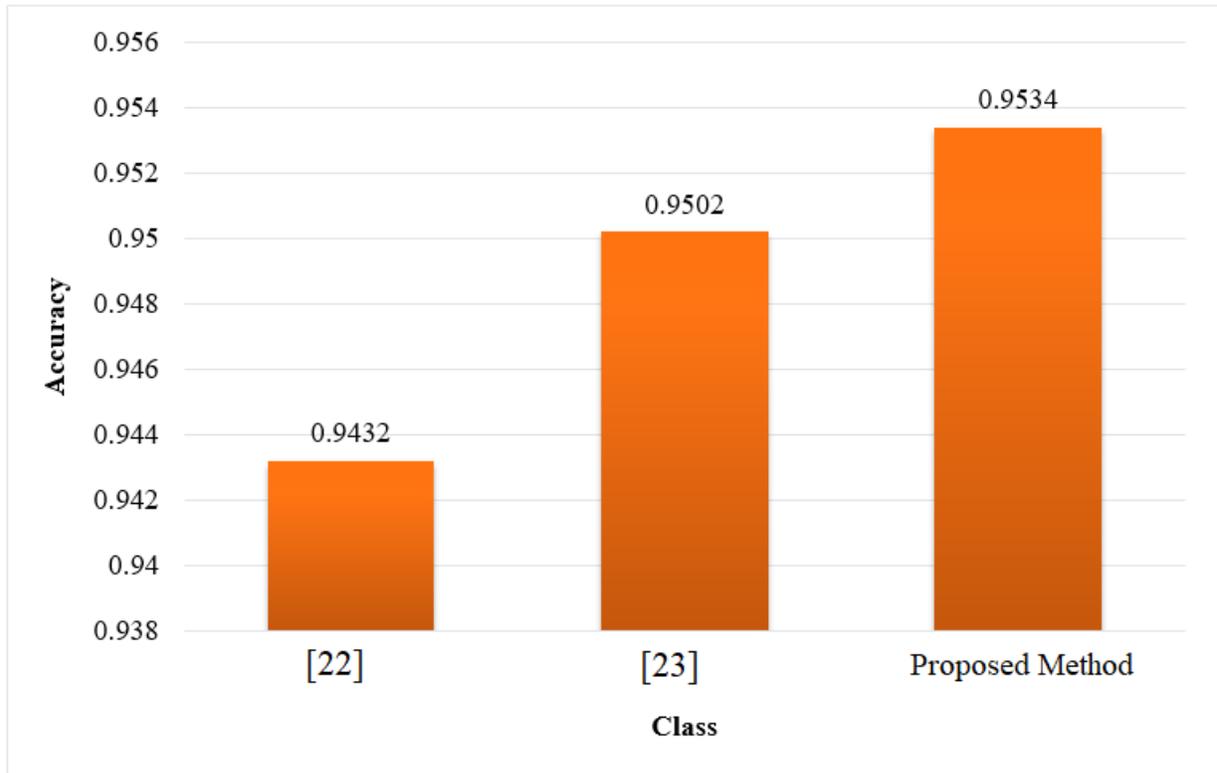

**Figure 11**: Comparison of the accuracy of the proposed method compared to other methods

Alt Text Figure 11[29 words]: Bar chart comparing the accuracy of three methods. Method [22] shows an accuracy of 0.9432, Method [23] has 0.9502, and the Proposed Method achieves the highest accuracy of 0.9534.

Figure 12 compares the True Positive Rate (TPR) of the proposed method with the two other methods. The proposed method consistently achieves a higher TPR, which means it is more effective at correctly identifying malicious activities in the dataset.

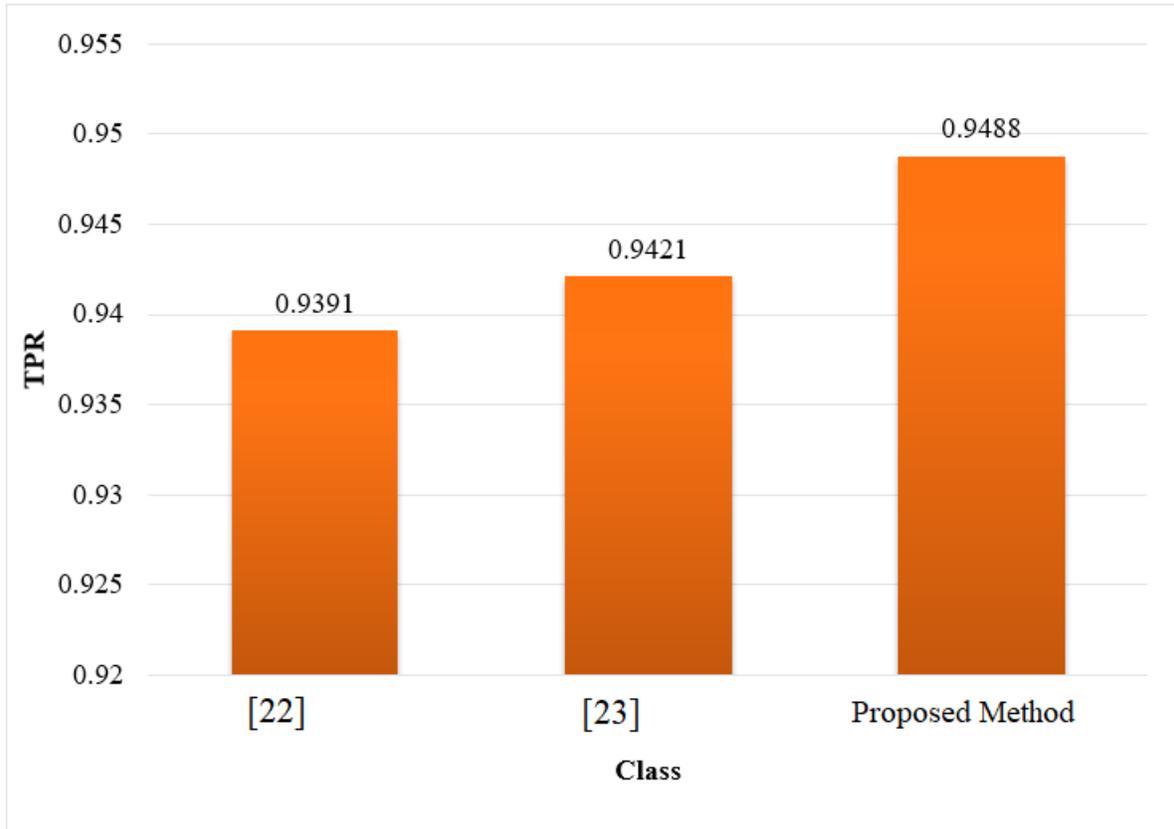

**Figure 12:** Comparison of the TPR rate of the proposed method and other methods

Alt Text Figure 12[35 words]: Bar chart comparing the TPR (True Positive Rate) of three methods. Method [22] has a TPR of 0.9391, Method [23] has a TPR of 0.9421, and the Proposed Method achieves the highest TPR of 0.9488

Figure 13 compares the False Positive Rate (FPR) of the proposed method with the other methods. A lower FPR means fewer normal (non-attack) activities are mistakenly flagged as attacks. As the graph shows, the proposed method has a lower average FPR than the other two methods, indicating that it is more accurate in avoiding false alarms.

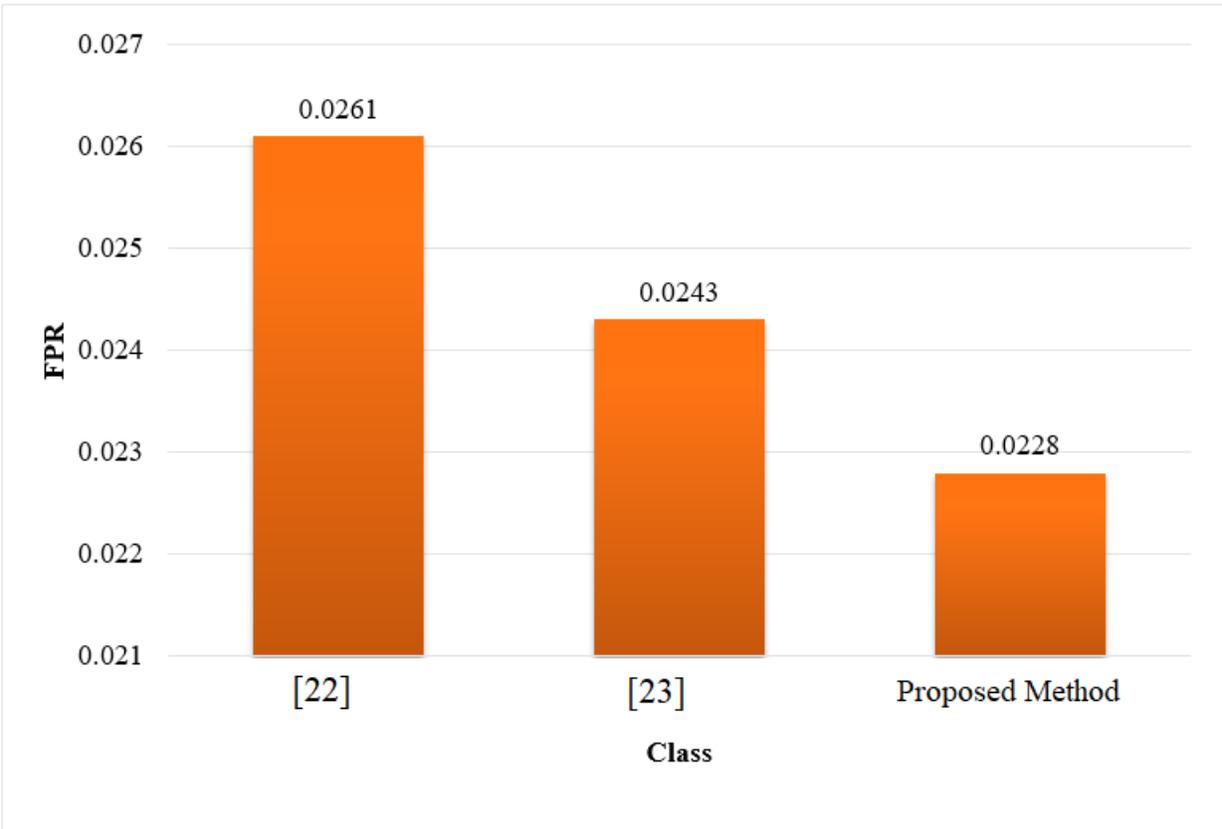

**Figure 13:** Comparison of the FRR rate of the proposed method and other methods

Alt Text Figure 13[36 words]: Bar chart comparing the FRR (False Rejection Rate) of three methods. Method [22] has an FRR of 0.0261, Method [23] has an FRR of 0.0243, and the Proposed Method has the lowest FRR of 0.0228.

## Conclusion

To protect IoT networks, Intrusion Detection Systems (IDS) are employed to detect unusual activities and alert users to potential threats. When an intrusion is identified, the IDS can initiate reconfiguration processes to isolate or remove compromised devices from the network, preserving overall network integrity. Proposed method aims to enhance security in IoT networks by implementing a more efficient IDS. The proposed method combines the Grasshopper Optimization Algorithm (GOA) for feature selection with Support Vector Machine (SVM) for attack classification. The system is tested using the NSL-KDD dataset, a refined version of the KDD99 dataset commonly used for evaluating intrusion detection models. After data gathering and pre-processing—which includes standardizing data and removing outliers—the system applies feature selection to reduce dimensionality and improve classifier accuracy. In the feature selection phase, GOA identifies the most critical features, minimizing resource consumption while maximizing detection accuracy. These selected features are then used by the SVM classifier to categorize network traffic as either normal or malicious. The proposed model demonstrates impressive performance in detecting intrusions, achieving superior results across key metrics, including accuracy, True Positive Rate (TPR), and False Positive Rate (FPR), compared to existing methods. By harnessing the optimization strengths of GOA alongside the robust classification capabilities of SVM, the system effectively identifies security threats while reducing false alarms.

**Authorship contribution statement**

Shiva Sattarpour: Conception and design of study/review/case series, Acquisition of data, laboratory search, Drafting of article.

Ali Barati: Conception and design of study/review/case series, Analysis and interpretation of data collected, Final approval and guarantor of manuscript.

Hamid Barati: Acquisition of data, Laboratory or clinical/literature search, Analysis and interpretation of data collected, Final approval.

Declaration of competing interest

The authors declare that there is no conflict of interests regarding the publication of this manuscript.

Data availability

No data was used for the research described in the article.

Acknowledgments

None.

Ethics approval

This article does not contain any studies with human participants.

Funding Declaration

No Funding